\documentclass[12pt]{article}
\usepackage{amsmath,amssymb}
\usepackage{latexsym}
\usepackage{graphicx}
\usepackage[labelfont=bf]{caption}
\newtheorem{thm}{Theorem}[section]

\newtheorem{defn}{Definition}[section]

\newtheorem{rem}{Remark}[section]
\newtheorem{rems}{Remarks}[section]

\newfont{\sss}{cmss12 scaled 1000}

\title{A recursive construction for simple $t-$designs using resolution}
\author{ Tran van Trung  \\
         Institut f\"ur Experimentelle Mathematik \\ 
         Universit\"at Duisburg-Essen \\
         Thea-Leymann-Stra\ss e 9, 45127 Essen, Germany }

\date{}

\begin{document}                  
                
\maketitle

\begin{abstract}
 This work presents a recursive construction for simple $t$-designs
 using resolutions of the ingredient designs. The result extends 
 a construction of $t$-designs in our recent paper 
 \cite{TvT2016}. Essentially, the method in \cite{TvT2016} describes 
 the blocks of a constructed design as a collection of block unions
 from a number of appropriate pairs of disjoint ingredient designs.
 Now, if some pairs of these ingredient $t$-designs 
 have both a suitable $s$-resolution, then we can define 
 a distance mapping  on their resolution classes. 
 Using this mapping enables us to have more possibilities for 
 forming blocks from those pairs. The method makes it possible
 for constructing many new simple $t$-designs. 
 We give some application results of the new construction.  
 
 \end{abstract}

\vspace{0.1in}\noindent
{\bf 2010 Mathematics Subject Classification:} 05B05

\vspace{0.1in}\noindent
{\bf Keywords:} recursive construction, resolution, large set, 
   simple $t-$design.


\section{Introduction}
  In a recent paper \cite{TvT2016} we have presented a recursive method for
  constructing simple $t$-designs for arbitrary $t$. 
  The method is of combinatorial nature since
  it requires finding solutions for the indices of the ingredient 
  designs that satisfy a certain set of equalities. In essence, the core 
  of the construction is that the blocks of a constructed design 
  are built as a collection of block unions from 
  a number of appropriate 
  pairs of disjoint ingredient designs. In particular, when a
  pair of ingredient designs is used, we take as new blocks the unions 
  of all the pairs of blocks in the two ingredient designs. For the sake of
  simplicity we refer to this construction method as the basic method
  or the basic construction. 
  
  In the present paper we describe an extension of the basic 
  construction by assuming
  that a subset of pairs of ingredient designs have
  suitable resolutions. For those given  
  pairs we may define a distance mapping on their resolution classes.  
  By using this mapping we have more possibilities for
  forming blocks from those pairs other than taking
  the unions of all possible pairs of blocks in the ingredient designs. 
  This construction actually extends the basic construction since
  many new simple $t$-designs can only be constructed with the 
  new method. The crucial point of this extension is the use 
  of $s$-resolutions for $t$-designs. The concept of $s$-resolution
  may be viewed as a generalization of the notion of parallelism,
  which may be termed as $(1,1)$-resolution,
  i.e. the blocks of the $t$-design can be partitioned into classes of
  mutually disjoint blocks such that every point is in exactly one
  block of each class. To date very little is known about
  $s$-resolutions for $t$-designs when $s \geq 2$, except for
  the trivial $t$-designs. In this case, an $s$-resolution
  of the trivial $t$-design
  turns out to be a large set of $s$-designs. A great deal
  of results about large sets of $s$-designs have been achieved
  by many researchers, see the references below.
  We will describe our construction in terms of $s$-resolution
  for $t$-designs in general. However we will restrict 
  its applications
  just for the case where pairs of trivial designs are used 
  and each has a suitable large set. Even with this
  limitation we find that the construction using resolution 
  still possesses its strength
  since many simple $t$-designs can be constructed. 
     
 It is worthwhile to emphasize that constructing simple $t-$designs for 
 large $t$ is a challenging problem in design theory.  
 There are several major approaches to the problem. These include
constructing $t-$designs from large sets of $t-$designs, for instance
\cite{Ajoodani96, Khos-Ajoodani91, Dehon83, Hartman87, Khos-Tayfeh2006,
 Kramer93, Kreher93, Laue2001, Laue2007, Teirlinck84, 
Teirlinck87, Teirlinck89, QRWu91};
constructing $t-$designs by using prescribed automorphism groups,
for example   \cite{Alltop69, Alltop72, betten95, betten98, betten99,
 bier95a, bier95b, Denniston76, Kramer76, Kreher90, Mag84, Schmalz93}; 
or contructing $t-$designs via recursive construction methods,
 see for instance
\cite{driessen78, Jimbo2011, Mag87, Tuan2002, 
Sebille2001, TvT84, TvT86, TvT2001, TvT2016, TvT2016B}.

\section{Preliminaries}
  
We recall some basic definitions. A $t-$design, denoted by
$t-(v,k,\lambda)$, is a pair $(X, {\mathfrak B})$, where $X$ is a 
$v-$set of {\it points} and ${\mathfrak B}$ is a collection of $k-$subsets,
called {\it blocks}, of $X$ having the property that every $t-$set 
of $X$ is a subset of exactly $\lambda$ blocks in ${\mathfrak B}$.
The parameter $\lambda$ is called the {\it index} of the design. 
A $t-$design is called {\it simple} if no two blocks are
 identical
i.e. no block of ${\cal B}$ is repeated; otherwise, it is called
{\it non-simple} (i.e. ${\mathfrak B}$ is a multiset). 
It can be shown by simple counting that a
$t-(v,k,\lambda)$ design is an $s-(v,k,\lambda_s)$ design for 
$0 \leq s \leq t$,
where $\lambda_s = \lambda{v-s \choose t-s}/{k-s \choose t-s}.$
Since $\lambda_s$ is an integer, necessary 
conditions for the parameters of a $t-$design are
${k-s \choose t-s} | \lambda{v-s \choose t-s}$, for $0 \leq s \leq t.$
For given $t,v$ and $k$, we denote by 
$\lambda_{\mbox{min}}(t,k,v)$, or 
$\lambda_{\mbox{min}}$ for short, the smallest
positive integer such that these conditions are satisfied for
all $0 \leq s \leq t.$ By complementing each block in $X$ of a
$t-(v,k,\lambda)$ design, we obtain a $t-(v,v-k, \lambda^*)$ design with
$ \lambda^* = \lambda{v-k \choose t}/{k \choose t}$, hence we shall 
assume that $k \leq v/2.$ The largest value for $\lambda$ for which
a simple $t-(v,k, \lambda)$ design exists is denoted 
by $\lambda_{\mbox{max}}$
and we have $\lambda_{\mbox{max}}={v-t \choose k-t}.$ The simple
$t-(v,k,\lambda_{\mbox{max}})$ design is called the {\it complete}
design or the {\it trivial} design. A $t-(v,k,1)$ design is called a
{\em t-Steiner system.}

We refer the reader to \cite{bjl99, handbook07} 
for more information about designs.

 \begin{defn}\label{def1}
  A $t-(v,k,\lambda)$-design $(X, {\mathfrak B})$ is
 said to be $(s,\tau)$-resolvable with $ 0 < s < t $, 
 if its block set $\mathfrak B$ 
 can be partitioned into $N$ classes 
 ${\mathfrak A}_1, \ldots, {\mathfrak A}_N$ such that
 $(X,{\mathfrak A}_i)$ is a $s-(v,k,\tau)$ design for all $i=1,\ldots, N.$
 Each ${\mathfrak A}_i$ is called a resolution class. We also say that 
 a $t-(v,k,\lambda)$-design has an $s-$resolution, 
 if it is $(s,\tau)$-resolvable for a certain $\tau$.
 
 \end{defn}
 It is worth noting that the concept of resolvability 
 (i.e. $(1,1)$-resolvability) for BIBD introduced
 by Bose in 1942 \cite{Bose1942} was generalized by
 Shrikhande and Raghavarao to $\sigma$-resolvability (i.e. 
 $(1, \tau)$-resolvability) for BIBD in 1963 \cite{SR1963}.
 A definition of $(s, \lambda)$-resolvability 
 for t-designs with $t \geq 3$  may be found in
 \cite{Baker1976}. In that paper Baker shows that the Steiner quadruple
 system $3-(4^m,4,1)$ constructed from an even dimensional
 affine space over the field of two elements has a $(2,1)-$resolution.
 Also, Teirlinck shows for example that there
 exists a 2-resolvable $3-(2p^n +2,4,1)$ design with 
 $p \in \{7,31,127 \}$,
 for any positive integer $n$, \cite{Teirlinck94}. To date,
 very little is known about $s-$resolution of 
 non-trivial $t-(v,k,\lambda)$
 designs for $ t \geq 3$ and $s \geq 2.$  

 When $(X, {\mathfrak B})$ is the trivial $t-(v,k,{v-t \choose k-t})$ 
 design,
 then an $(s,\tau)-$resolution of $(X,{\mathfrak B})$ is called a 
 {\it large set}. Thus, a large set is a partition of the complete
 $t-(v,k,{v-t \choose k-t})$ design into $s-(v,k,\tau)$ designs, 
 and is denoted
 by LS$[N](s,k,v)$, where $N={v-s \choose k-s}/\tau$ is the number
 of resolution classes in the partition. 

 We define a distance on the resolution classes of a $t-$design
 as follows.

 \begin{defn} \label{distance}
  Let $D$ be a $t-(v,k,\lambda)$ design admitting
  an $(s,\tau)$-resolution with
  ${\mathfrak A}_1, \ldots , {\mathfrak A}_N$ as resolution classes.
  Define a distance  between any two classes 
  ${\mathfrak A}_i$ and ${\mathfrak A}_j$ by
  $d({\mathfrak A}_i, {\mathfrak A}_j)= \min \{|i-j|, N-|i-j| \}.$  

 \end{defn}


\subsection{The basic construction}

 In this section, we summarize the basic construction
 as described in \cite{TvT2016}. This preparation is
 necessary for the description of the construction using resolution
 in the next section. 

 \vspace{2mm}
 We first give notation and definitions.
 Let $t, \; v, \; k$ be non-negative integers such that 
 $v \geq k \geq t \geq 0$.
 Let $X$ be a $v$-set and let $X= X_1 \cup X_2$ be a partition of $X$
 (i.e $X_1 \cap X_2 = \emptyset$) with $|X_1|=v_1$ and $|X_2|=v_2$.  
 
 The parameter set
 $t-(v_2,j, \bar{\lambda}^{(j)}_t)$ for a design indicates that the point set 
 of the design is $X_2$. Also, a design defined on the point set $X_2$ 
 is denoted by $\bar{D}=(X_2, \bar{\mathfrak B})$.

 \begin{itemize}
 \item[(i)]
  For $i=0, \ldots, t$, let $D_i =(X_1, {\mathfrak B}^{(i)})$ 
   be the 
  complete $i-(v_1,i,1)$ design. For $i=t+1, \ldots, k$, let 
  $D_i =(X_1, {\mathfrak B}^{(i)})$ be a simple 
  $t-(v_1,i, \lambda^{(i)}_t)$ design.
 \item[(ii)]
  Similarly, for $i=0, \ldots, t$, let 
  $\bar{D}_i =(X_2, \bar{\mathfrak B}^{(i)})$ 
  be the  complete $i-(v_2,i,1)$ design. And for $i=t+1, \ldots, k$, 
  let   $\bar{D}_i =(X_2, \bar{\mathfrak B}^{(i)})$ be a 
  simple $t-(v_2,i, \bar{\lambda}^{(i)}_t)$ design.
 \item[(iii)]
  Two degenerate cases for designs occur when either $v=k=t=0$ or $v=k$. 
  The first case $v=k=t=0$ gives an ``empty''
  design, denoted by $\emptyset$, however we use the convention that the
  number of blocks of the empty design is 1 (i.e. the unique block is
  the empty block). The second case $v=k$ gives a degenerate $k$-design
  having just 1 block consisting of all $v$ points. 
  Thus, in these two extreme
  cases the number of blocks of the designs is always 1.

 \item[(iv)]
  We denote by $T_{(r,t-r)}$ a $t$-subset $T$ of $X$ with
  $|T\cap X_1|=r$ and hence $|T\cap X_2|=t-r$, for $r=0,\ldots, t$.
  It is clear that any $t$-subset of $X$ is a $T_{(r,t-r)}$ set 
  for some $r \in \{0, \ldots, t\}$.
 
  \item[(v)]
  Let $X$ be a finite set and let $u \in\{0,1\}$. The notation $X \times [u]$ 
  has the following meaning. $X \times [0]$ is  
  the empty set $\emptyset$, and  $X \times [1] = X$.

 \end{itemize}
      
 The basic construction in \cite{TvT2016} is as follows.

 Consider $(k+1)$ pairs of simple designs $(D_i,\bar{D}_{k-i})$
 for $i=0, \ldots, k$, where   
  $D_i =(X_1, {\mathfrak B}^{(i)})$ is a simple $t-(v_1,i, \lambda^{(i)}_t)$
  design and 
  $\bar{D}_{k-i} =(X_2, \bar{\mathfrak B}^{(k-i)})$
  a simple $t-(v_2,k-i, \bar{\lambda}^{(k-i)}_t)$ design,
  as defined above. 
 For each pair $(D_i,\bar{D}_{k-i})$ define
  $${\mathfrak B}_{(i,{k-i})} := \{ B= B_i \cup \bar{B}_{k-i} \; / \;
   B_i \in { \mathfrak B}^{(i)}, \bar{B}_{k-i} \in 
      \bar{\mathfrak B}^{(k-i)} \}. $$ 
 
 Define
 $${\mathfrak B} := {\mathfrak B}_{(0,k)}\times [u_0] \cup 
             {\mathfrak B}_{(1,k-1)}\times [u_1] 
             \cup   \cdots \cup
             {\mathfrak B}_{(k-1,1)} \times [u_{k-1}] 
             \cup {\mathfrak B}_{(k,0)}\times [u_k], $$ 
where $u_i \in \{0,1\}$, for $i=0, \ldots, k.$  

 It should be remarked that the notation 
 ${\mathfrak B}_{(i,k-i)}\times [u_i]$, as defined in (v) above, 
 indicates that either we
 have an empty set $\emptyset$ (when $u_i=0$) or the set 
 ${\mathfrak B}_{(i,k-i)}$ itself (when $u_i=1$).
 The empty set case means that the pair $(D_i,\bar{D}_{k-i})$
 is not used and the other case means the pair $(D_i,\bar{D}_{k-i})$
 is used.

 It can be shown that a given $t$-set $T_{(r,t-r)}$ of $X$ the number 
 of blocks in $\mathfrak B$ containing $T_{(r,t-r)}$ is equal to

\begin{eqnarray*}
 L_{r,t-r} & := & \sum_{i=0}^k u_i.\lambda^{(i)}_r. 
                        \bar{\lambda}^{(k-i)}_{t-r}. 
\end{eqnarray*} 
 
Therefore, if
\begin{eqnarray*}
 & & L_{0,t} = L_{1,t}=L_{2,t-2}= \cdots = L_{t,0}:=\Lambda, 
\end{eqnarray*}
 where $\Lambda$ is a positive integer, then $(X, {\mathfrak B})$
 forms a simple $t$-design with parameters
 $t-(v,k, \Lambda)$.

 We record the basic construction in following theorem.

\begin{thm}[Basic construction]\label{basic}
 Let v, k, t be integers with $v > k > t \geq 2$. Let $X$ be a $v$-set
 and let $X=X_1 \cup X_2$ be a partition of $X$ with $|X_1|=v_1$ and
 $|X_2|=v_2$. 
  Let $D_i =(X_1, {\mathfrak B}^{(i)})$ be the 
  complete $i-(v_1,i,1)$ design for $i=0, \ldots, t$  
  and let 
  $D_i =(X_1, {\mathfrak B}^{(i)})$ be a simple $t-(v_1,i, \lambda^{(i)}_t)$
  design for $i=t+1, \ldots, k$.
  Similarly, let 
  $\bar{D}_i =(X_2, \bar{\mathfrak B}^{(i)})$ 
  be the  complete $i-(v_2,i,1)$ design for $i=0, \ldots, t$, and let 
  $\bar{D}_i =(X_2, \bar{\mathfrak B}^{(i)})$ be a 
  simple $t-(v_2,i, \bar{\lambda}^{(i)}_t)$ design for $i=t+1, \ldots, k$.
  Define 
 $${\mathfrak B} = {\mathfrak B}_{(0,k)}\times [u_0] \cup 
             {\mathfrak B}_{(1,k-1)}\times [u_1] 
             \cup   \cdots \cup
             {\mathfrak B}_{(k-1,1)} \times [u_{k-1}] 
             \cup {\mathfrak B}_{(k,0)}\times [u_k], $$   
  where  
$${\mathfrak B}_{(i,{k-i})}= \{ B= B_i \cup \bar{B}_{k-i} \; / \;
   B_i \in {\mathfrak B}^{(i)}, \bar{B}_{k-i} \in 
     \bar{\mathfrak B}^{(k-i)} \}. $$ 
  Assume that 
\begin{eqnarray}\label{Equalities}
  & & L_{0,t}= L_{1,t-1}=L_{2,t-2}= \cdots = L_{t,0}:=\Lambda, 
\end{eqnarray}
  for a positive integer $\Lambda$, where  
  \begin{eqnarray}\label{L-values}
  L_{r,t-r} & = & \sum_{i=0}^k u_i.\lambda^{(i)}_r. 
    \bar{\lambda}^{(k-i)}_{t-r},
  \end{eqnarray}
  $r=0, \ldots , t$, and $u_i\in\{0,1\}$, for $i=0, \ldots, k.$ 
 Then $(X, {\mathfrak B})$ is a simple $t-(v,k, \Lambda)$ design. 
\end{thm}

\section{The construction using resolution}

In this section we describe a recursive construction
of simple $t-$designs using resolution. Note that in the 
basic construction, if a pair 
$(D_i, \bar{D}_{k-i})$ is used in the construction  
(i.e. $u_i=1$), then the new blocks formed by this pair
consist of taking the union of each block of $D_i$
with each block of $\bar{D}_{k-i}$. 
The crucial idea of the construction using resolution is that
if $D_i$ and $\bar{D}_{k-i}$  
have appropriate $s_1-$ and $s_2-$resolutions with the same 
number of resolution classes, then the new blocks are formed
according to the distance mapping defined on the resolution
classes of $D_i$ and $\bar{D}_{k-i}$ rather than taking
the unions of each block of $D_i$ with each block of $\bar{D}_{k-i}$.

In the following we go into detail of the construction.
We make use of the notation and definitions for the basic construction
in the previous section.  
When for a certain $i\in \{0, \ldots, k\}$
the $t-(v_1,i, \lambda_t^{(i)})$ design  $D_i=(X_1, {\mathfrak B}^{(i)})$
has an $s_i-$resolution, i.e. $D_i$ can be partitioned into $N_i$
disjoint $(X_1,{\mathfrak A}^{(i)}_h)$ designs with parameters
$s_i-(v_1, i, \lambda_{s_i}^{*(i)})$, $s_i < t$, then we write
$$ {\mathfrak B}^{(i)}= \bigcup_{h=1}^{N_i} {\mathfrak A}^{(i)}_h, $$
where
$$ N_i ={\lambda_t^{(i)}} {v_1-s_i \choose t-s_i}/
         {\lambda_{s_i}^{*(i)}} {i-s_i \choose t-s_i}. $$

Similarly, we  write
$$ \bar{\mathfrak B}^{(k-i)}= 
   \bigcup_{h=1}^{{\bar N}_{k-i}} \bar{\mathfrak A}^{(i)}_h, $$
when the blocks of a
$t-(v_2,k-i, {\bar\lambda_t}^{(k-i)})$ design  
 ${\bar D}_{k-i}=(X_2, \bar{\mathfrak B}^{(k-i)})$  
can  be partioned into ${\bar N}_{k-i}$ disjoint
$(X_2,\bar{\mathfrak A}^{(k-i)}_h)$ designs with parameters 
$s_{k-i}-(v_2, k-i, \bar\lambda_{s_{k-i}}^{*(k-i)})$, where
$$ {\bar N}_{k-i} ={\bar\lambda}_t^{(k-i)} {v_2-s_{k-i} \choose t-s_{k-i}}/
   {\bar\lambda}_{s_{k-i}}^{*(k-i)} {k-i-s_{k-i} \choose t-s_{k-i}} $$ 
is the number of $s_{k-i}-$resolution classes.

Let $ K= \{(0,k),(1, k-1), \ldots, (k-1,1), (k,0) \}.$
Assume there exists a subset $R \subseteq K$ such that if $(i,k-i) \in R$,
then $D_i$ and ${\bar D}_{k-i}$ have $s_i-$resolution of size $N_i$
and $s_{k-i}-$resolution of size ${\bar N}_{k-i}$, respectively, 
satisfying the following conditions.

\begin{itemize}
\item[(i)] $N_i = {\bar N}_{k-i},$

\item[(ii)]
    $ s_i + s_{k-i} \geq 2\lfloor \frac{t}{2} \rfloor.$ 

\end{itemize}

The construction of $t-$designs  using resolution consists 
of building two types of blocks.

\begin{itemize}
\item [(1)]
   
   For each pair $(i, k-i) \in K \setminus R$ form a subset
   of new blocks ${\mathfrak B}_{(i,{k-i})}$ 
   from the pair $(D_i, {\bar D}_{k-i})$ as
   $${\mathfrak B}_{(i,{k-i})}:= \{ B= B_i \cup \bar{B}_{k-i} \; / \;
   B_i \in {\mathfrak B}^{(i)}, \bar{B}_{k-i} \in 
     \bar{\mathfrak B}^{(k-i)} \}. $$ 
 
 \item [(2)] 

   For each pair $(i, k-i) \in R$ form a subset of new blocks
   ${\mathfrak B}^*_{(i,{k-i})}$ from $(D_i, {\bar D}_{k-i})$
   by using $s_i-$resolution of $D_i$ and $s_{k-i}-$resolution
   of ${\bar D}_{k-i}$ as follows. 
 $${\mathfrak B}^*_{(i,{k-i})}:= \{ B_i \cup {\bar B}_{k-i} /  B_i \in 
    {\mathfrak A}^{(i)}_h, {\bar B}_{k-i} \in \bar{\mathfrak A}^{(k-i)}_j,\; 
        \varepsilon_i \leq 
    d({\mathfrak A}^{(i)}_h,{\mathfrak A}^{(i)}_j) \leq w_i,\varepsilon_i =0,1;
      \;  w_i \leq \lfloor \frac{N_i}{2} \rfloor  \}. $$
\end{itemize}

Further, define
 $$z_i :=(2w_i+1-\varepsilon_i), \mbox{ if }  \; w_i < \frac{N_i}{2},
 \mbox{ and }
 z_i :=(2w_i-\varepsilon_i), \mbox{ if }  \; w_i = \frac{N_i}{2}. $$
 Note that $w_i$ and  $z_i$ are considered as variables.
 
Now, let $T_{(r, t-r)}$ be a $t-$set of $X$ for $r=0, \ldots , t$.
According to the property of $s_i$ and $s_{k-i}$ one of the following
cases has to occur.

\begin{itemize}
\item [(a)] 
$r \leq s_i$ and $t-r \leq s_{k-i}$. Then $T_{(r, t-r)}$
 is contained in
$$\Lambda^{*(i,k-i)}_{r,t-r} = \lambda^{*(i)}_r. 
              {\bar\lambda}^{*(k-i)}_{t-r}. N_i. z_i$$
blocks of ${\mathfrak B}^*_{(i,{k-i})}.$
\item [(b)]
 $r \leq s_i$ and $t-r > s_{k-i}$. Then $T_{(r, t-r)}$
 is contained in
$$\Lambda^{*(i,k-i)}_{r,t-r} = \lambda^{*(i)}_r. {\bar\lambda}^{(k-i)}_{t-r}.
                   z_i$$
blocks of ${\mathfrak B}^*_{(i,{k-i})}.$

\item [(c)]
$r > s_i$ and $t-r \leq  s_{k-i}$. Then $T_{(r, t-r)}$
 is contained in
$$\Lambda^{*(i,k-i)}_{r,t-r} = \lambda^{(i)}_r. {\bar\lambda}^{*(k-i)}_{t-r}.
                   z_i$$
blocks of ${\mathfrak B}^*_{(i,{k-i})}.$
\end{itemize}

\vspace{4mm}
 It is straightforward to verify the values of 
 $\Lambda^{*(i,k-i)}_{r,t-r}$ for the cases (a), (b) and (c) above.
 In case (a) we have that each $r-$subset of $X_1$ is contained in 
 $\lambda^{*(i)}_r$ blocks of $ {\mathfrak A}^{(i)}_h$ and each 
 $(t-r)-$subset of $X_2$  is contained in 
 ${\bar\lambda}^{*(k-i)}_{t-r}$ blocks of $\bar{\mathfrak A}^{(k-i)}_j$.
 Thus each pair
 $({\mathfrak A}^{(i)}_h, \bar{\mathfrak A}^{(k-i)}_j)$ contributes 
 $\lambda^{*(i)}_r.{\bar\lambda}^{*(k-i)}_{t-r}$ new blocks to
 ${\mathfrak B}^*_{(i,{k-i})}$. Now each of the $N_i$ resolution classes 
 ${\mathfrak A}^{(i)}_1, \ldots, {\mathfrak A}^{(i)}_{N_i}$ 
 is combined with $z_i$ resolution classes of
 $\bar{\mathfrak A}^{(k-i)}_1, \ldots, \bar{\mathfrak A}^{(k-i)}_{N_i}$, 
 therefore we have    
 $\Lambda^{*(i,k-i)}_{r,t-r} = \lambda^{*(i)}_r. 
              {\bar\lambda}^{*(k-i)}_{t-r}. N_i. z_i$

  In case (b) each $r-$subset of $X_1$ is contained in
 $\lambda^{*(i)}_r$ blocks of $ {\mathfrak A}^{(i)}_h$ and each
 $(t-r)-$subset of $X_2$  is contained in exactly
 ${\bar\lambda}^{(k-i)}_{t-r}$ blocks of $\bar{\mathfrak B}^{(k-i)}$.
 These blocks are distributed in the $N_i$ resolution classes
 $\bar{\mathfrak A}^{(k-i)}_1, \ldots, \bar{\mathfrak A}^{(k-i)}_{N_i}$. 
 Each class $\bar{\mathfrak A}^{(k-i)}_j$ is combined
 $z_i$ times with ${\mathfrak A}^{(i)}_h$. Hence, in this case,
 the contribution of the new blocks to ${\mathfrak B}^*_{(i,{k-i})}$  is 
$\Lambda^{*(i,k-i)}_{r,t-r} = \lambda^{*(i)}_r. {\bar\lambda}^{(k-i)}_{t-r}.
                   z_i.$
 
\vspace{2mm} 
 The case (c) is similar to case (b).

\vspace{2mm}
Define
$$ \mathfrak B :=\bigcup_{(i,k-i)\in R} {\mathfrak B}^*_{(i,{k-i})}\times [u_i]
\bigcup_{(i,k-i)\in k\setminus R} {\mathfrak B}_{(i,{k-i})}\times [u_i],$$
 with $u_i\in \{0,1\}$, $i =0, \ldots , k.$

\vspace{2mm}
The above presentation can be summarized as follows.
Let $T_{(r,t-r)}$  be a $t-$subset of $X$ for $r=0, \ldots, t$. 
The entire number of new blocks in ${\mathfrak B}_{(i,k-i)}$
containing $T_{(r,t-r)}$, for all $(i,k-i) \in K \setminus R$, is then
\begin{eqnarray*}
 \sum_{(i,k-i) \in K \setminus R} u_i.\Lambda_{(r,t-r)}^{(i,k-i)} 
 & = & \sum_{(i,k-i) \in K \setminus R} u_i. \lambda_r^{(i)}. 
                     \bar{\lambda}_{t-r}^{(k-i)}.
\end{eqnarray*}

The entire number of new blocks in ${\mathfrak B}_{(i,k-i)}^*$
containing $T_{(r,t-r)}$, for all $(i,k-i) \in R$, is then
 $$\sum_{(i,k-i) \in R} u_i.\Lambda_{(r,t-r)}^{*(i,k-i)}, $$
where
$$\Lambda_{(r,t-r)}^{*(i,k-i)} = \left\{ 
   \begin{array}{rl}
    \lambda^{*(i)}_r.{\bar\lambda}^{*(k-i)}_{t-r}. N_i. z_i &
      \text{if } r \leq s_i, \; t-r\leq s_{k-i}, \\
      \lambda^{*(i)}_r. {\bar\lambda}^{(k-i)}_{t-r}.z_i  &
      \text{if } r \leq s_i, \; t-r > s_{k-i}, \\
      \lambda^{(i)}_r. {\bar\lambda}^{*(k-i)}_{t-r}. z_i &
      \text{if } r > s_i, \; t-r \leq s_{k-i}. 
  \end{array} \right.
$$
It follows that the number of blocks in $\mathfrak B$ containing
$T_{r,t-r}$ is equal to
$$ L_{r, t-r} := 
  \sum_{(i,k-i) \in R} u_i.\Lambda_{(r,t-r)}^{*(i,k-i)} +
 \sum_{(i,k-i) \in K \setminus R} u_i.\Lambda_{(r,t-r)}^{(i,k-i)}. $$

Since any $t-$subset of $X$ is of form $T_{r,t-r}$ for some
$r \in \{ 0, \ldots, t  \}$, we see that
if 
$$L_{0,t} = L_{1,t-1} = \dots = L_{t,0} := \Lambda $$
for a positive integer $\Lambda$, then $(X, \mathfrak B)$ forms
a simple $t-$design with parameters $t-(v, k, \Lambda).$

We record the construction above in the following theorem.

\begin{thm}\label{main}
 Let v, k, t be integers with $v > k > t \geq 2$. Let $X$ be a $v$-set
 and let $X=X_1 \cup X_2$ be a partition of $X$ with $|X_1|=v_1$ and
 $|X_2|=v_2$. 
  Let $D_i =(X_1, {\mathfrak B}^{(i)})$ be the 
  complete $i-(v_1,i,1)$ design for $i=0, \ldots, t$  
  and let 
  $D_i =(X_1, {\mathfrak B}^{(i)})$ be a simple $t-(v_1,i, \lambda^{(i)}_t)$
  design for $i=t+1, \ldots, k$.
  Similarly, let 
  $\bar{D}_i =(X_2, \bar{\mathfrak B}^{(i)})$ 
  be the  complete $i-(v_2,i,1)$ design for $i=0, \ldots, t$, and let 
  $\bar{D}_i =(X_2, \bar{\mathfrak B}^{(i)})$ be a 
  simple $t-(v_2,i, \bar{\lambda}^{(i)}_t)$ design for $i=t+1, \ldots, k$.
  Let $K =\{(0,k), (1,k-1), \dots, (k-1,1), (k,0) \}.$ Suppose there
  exists a subset $R \subseteq K$ such that for each $(i, k-i) \in R$,
  the designs $D_i$ and $\bar D_{k-i}$ have $s_i-$resolution with $N_i$
  classes and $s_{k-i}-$resolution with $\bar N_{k-i}$ classes, respectively,
  satisfying the following conditions.
  \begin{itemize}
   \item [(i)]
    $N_i = \bar N_{k-i},$
    \item [(ii)]
     $s_i + s_{k-i} \geq 2 \lfloor \frac{t}{2} \rfloor.$
  \end{itemize}
  Define
$$ \mathfrak B =\bigcup_{(i,k-i)\in R} {\mathfrak B}^*_{(i,{k-i})}\times [u_i]
\bigcup_{(i,k-i)\in k\setminus R} {\mathfrak B}_{(i,{k-i})}\times [u_i],$$
for $u_i \in \{0,1\}$, $i=0, \ldots, k$,
 $${\mathfrak B}^*_{(i,{k-i})}:= \{ B_i \cup {\bar B}_{k-i} /  B_i \in 
    {\mathfrak A}^{(i)}_h, {\bar B}_{k-i} \in \bar{\mathfrak A}^{(k-i)}_j,\; 
        \varepsilon_i \leq 
      d({\mathfrak A}^{(i)}_h,{\mathfrak A}^{(i)}_j) \leq w_i,  
       \varepsilon_i =0,1;
      \;  w_i \leq \lfloor \frac{N_i}{2} \rfloor  \}, $$
with $w_i$ as variable,
where  $\mathfrak A^{(i)}_1, \ldots, \mathfrak A^{(i)}_{N_i}$ are
$s_i-$resolution classes of $D_i$, with $(X_1,\mathfrak A^{(i)}_h)$ as 
an $s_i-(v_1,i, \lambda^{*(i)}_{s_i})$ design;  and
 $\bar{\mathfrak A}^{(k-i)}_1, \ldots, \bar{\mathfrak A}^{(k-i)}_{N_i}$
 are  $s_{k-i}-$resolution classes of $\bar D_{k-i}$, with
 $(X_2,\bar{\mathfrak A}^{(k-i)}_h)$ as an  
 $s_{k-i}-(v_2,k-i, \bar{\lambda}^{*(k-i)}_{s_{k-i}})$ design; 
and
   $${\mathfrak B}_{(i,{k-i})}:= \{ B= B_i \cup \bar{B}_{k-i} \; / \;
   B_i \in {\mathfrak B}^{(i)}, \bar{B}_{k-i} \in 
     \bar{\mathfrak B}^{(k-i)} \}. $$ 
Define
$$ L_{r, t-r} := 
  \sum_{(i,k-i) \in R} u_i.\Lambda_{(r,t-r)}^{*(i,k-i)} +
 \sum_{(i,k-i) \in K \setminus R} u_i.\Lambda_{(r,t-r)}^{(i,k-i)}, $$
 for $r=0, \ldots, t $, where 
$$\Lambda_{(r,t-r)}^{*(i,k-i)} = \left\{ 
   \begin{array}{rl}
    \lambda^{*(i)}_r.{\bar\lambda}^{*(k-i)}_{t-r}. N_i. z_i &
      \text{if } r \leq s_i, \; t-r\leq s_{k-i}, \\
      \lambda^{*(i)}_r. {\bar\lambda}^{(k-i)}_{t-r}.z_i  &
      \text{if } r \leq s_i, \;  t-r > s_{k-i}, \\
      \lambda^{(i)}_r. {\bar\lambda}^{*(k-i)}_{t-r}. z_i &
      \text{if } r > s_i, \;  t-r \leq s_{k-i}. 
  \end{array} \right.
$$
with
 $z_i =(2w_i+1-\varepsilon_i), \mbox{ if }  \; w_i < \frac{N_i}{2},
 \mbox{ and }
 z_i =(2w_i-\varepsilon_i), \mbox{ if }  \; w_i = \frac{N_i}{2}; $
 and
 $$ \Lambda_{(r,t-r)}^{(i,k-i)} 
       = \lambda_r^{(i)}. \bar{\lambda}_{t-r}^{(k-i)}.$$
Assume that
\begin{eqnarray}\label{Equalities2} 
 & & L_{0,t} = L_{1,t-1} = \dots = L_{t,0} := \Lambda
\end{eqnarray}
for a positive integer $\Lambda$, then $(X, \mathfrak B)$ is 
a simple $t-(v, k, \Lambda)$ design.

\end{thm}
 
\begin{rems}
{\rm  
 \begin{enumerate}
 \item
 In the basic construction the set 
 $\mathfrak B_{(i,k-i)}$ of the new blocks 
 is uniquely determined as the unions of all the pairs of blocks
 in $D_i$ and  $\bar D_{k-i}.$
 Whereas in the construction using resolution in Theorem \ref{main}
 the set $\mathfrak B^*_{(i,k-i)}$ is no longer unique. 
 Its size varies according to the variable $z_i$.
 \item
 Theorem \ref{main} does not restrict to constructing simple
 $t-$ designs. Obviously, if any of the ingredient designs is non-simple,
 then the construction will yield non-simple designs.

\end{enumerate}
 
}
\end{rems}

\section{Applications}
 In this section we illustrate the construction 
 in Theorem \ref{main} through a number of 
 examples which show the strength of the method.

   In the following we will employ the notation from
   Chapter 4 : $t$-Designs with $t \geq 3$ 
   of the Handbook of Combinatorial Designs.  
   The parameter set $t-(v,k, \lambda)$ of a design will be written as
   $t-(v,k, m\lambda_{\text{min}}).$ 
   Since the supplement of a simple
   $t-(v,k, \lambda)$ design is a $t-(v,k, \lambda_{\text{max}}-\lambda)$
   design, we usually consider simple $t-(v,k, \lambda)$ designs with
   $\lambda \leq \lambda_{\text{max}}/2$.  Thus, the upper limit
   of $m$ of a constructed design will be 
$\text{LIM}= \lfloor \lambda_{\text{max}}/(2\lambda_{\text{min}}) \rfloor.$ 
But, it should be remarked that, when an ingredient design with 
   index $\lambda$ is used, then $\lambda$ can take on
   all possible values, i.e. 
   $\lambda_{\text{min}} \leq  \lambda \leq \lambda_{\text{max}}$.
  
\subsection{ Simple $5-(38, k, \Lambda)$ designs with $k=8,9,10$}
  We apply the construction in Theorem \ref{main} to the cases
  $t=5$, $v_1=v_2=19$ and $k=8,9,10.$

\subsubsection{Simple $5-(38,8, \Lambda)$ designs}
Here we show a detailed example to illustrate the construction.
 
 Let $ X= X_1 \cup X_2$ be a partition of the point set $X$ with 
 $|X|=38$ into two subsets $X_1$ and $X_2$ with $|X_1|=|X_2|=19.$ 
 For $i=0,1,2,3,4,5$ let $D_i=(X_1, {\cal B}^{(i)})$ be the complete
 $i-(19,i,\lambda^{(i)}_i):= i-(19,i,1)$ design. For $i=6,7,8$ let 
 $D_i=(X_1, {\cal B}^{(i)})$ be a simple $5-(19,i, \lambda^{(i)}_5)$
 design. These designs have the following parameters.   
   \begin{itemize}
    \item[$\bullet$]
     $5-(19,6,\lambda^{(6)}_5) = 5-(19,6,m2)$, $m=1,2,\ldots, 7$. 
    \item[$\bullet$]
    $5-(19,7,\lambda^{(7)}_5) =  5-(19,7,m7)$, $m=1,2,\ldots, 13$ 
    \item[$\bullet$]
    $5-(19,8,\lambda^{(8)}_5) =  5-(19,8,m28)$, $m=1,2,\ldots, 13$
   \end{itemize}
  Correspondingly, let
  ${\bar D_i}=(X_2, {\bar{\cal B}}^{(i)})$ be simple
  designs defined on $X_2$. 
  Here 
 $K=\{(0,8), (1,7), (2,6), (3,5), (4,4), (5,3), (6,2), (7,1), (8,0) \}.$
  
  It is known that the complete designs $D_i$ and $\bar{D}_i$
  for $i=3,4,5$ have each a 2-resolution with the number of resolution
  classes $N_i=17$, i.e. the large sets $\text{LS}[17](2,i,19)$,
  see for instance \cite{handbook07}. We choose
   $$R=\{(3,5), (4,4), (5,3) \}.$$ 
  Thus we have
 \begin{itemize}
\item
  ${\mathfrak B}^{(3)} =\bigcup_{j=1}^{17} {\mathfrak A}^{(3)}_j,$ 
   where $(X_1,{\mathfrak A}^{(3)}_j)$ is
  a $2-(19,3, \lambda^{*(3)}_2) = 2-(19,3,1)$ design, and
  $\lambda^{(3)}_3=1, \; 
   \lambda^{*(3)}_2=1, \; \lambda^{*(3)}_1=9, \; \lambda^{*(3)}_0=57;$
 \item  
  ${\mathfrak B}^{(4)} =\bigcup_{j=1}^{17} {\mathfrak A}^{(4)}_j,$ 
    where $(X_1,{\mathfrak A}^{(4)}_j)$ is
  a $2-(19,4, \lambda^{*(4)}_2) = 2-(19,4,8)$ design and
  $\lambda^{(4)}_4=1, \; \lambda^{(4)}_3=16, \; 
  \lambda^{*(4)}_2=8, \; \lambda^{*(4)}_1=48, \; \lambda^{*(4)}_0=228;$
\item
  ${\mathfrak B}^{(5)} =\bigcup_{j=1}^{17} {\mathfrak A}^{(5)}_j,$ 
  where $(X_1,{\mathfrak A}^{(i)}_j)$ is
  a $2-(19,5, \lambda^{*(4)}_2) = 2-(19,5,40)$ design and
  $\lambda^{(5)}_5=1, \; \lambda^{(5)}_4=15, \; \lambda^{(5)}_3=120, \; 
   \lambda^{*(5)}_2=40, \; \lambda^{*(5)}_1=180, \; \lambda^{*(5)}_0=684;$

\end{itemize}
  Similarly, the complete designs $\bar{D}_i$ have the same 2-resolution
  as $D_i$, each having $\bar{N}_i=17$ resolution classes,
   for $i=3,4,5.$ Thus
  $\bar{\mathfrak B}^{(i)} = \bigcup_{j=1}^{17} \bar{\mathfrak A}^{(i)}_j,$ 
  and each
  $(X_2,\bar{\mathfrak A}^{(i)}_j)$ 
  is a $2-(19,i, \bar{\lambda}^{*(i)}_2)$ design
  with $\bar{\lambda}^{*(i)}_2 = \lambda^{*(i)}_2.$
   
  We compute 
  \begin{eqnarray*} 
   L_{r, 5-r}
      & = & \sum_{(i,8-i) \in R} u_i.\Lambda_{(r,5-r)}^{*(i,8-i)} +
 \sum_{(i,8-i) \in K \setminus R} u_i.\Lambda_{(r,5-r)}^{(i,8-i)},
  \end{eqnarray*}
 for $r=0, \ldots, 5 $,  and $u_i =0,1.$ 
  If $(i, 8-i) \in K \setminus R$, then 
  $$\Lambda_{(r,5-r)}^{(i,8-i)}=\lambda^{(i)}_r.\bar{\lambda}^{(8-i)}_{5-r}.$$
  If $(i, 8-i) \in R$, then the values of $\Lambda_{(r,5-r)}^{*(i,8-i)}$
  are computed by using the formula 
$$\Lambda_{(r,t-r)}^{*(i,k-i)} = \left\{ 
   \begin{array}{rl}
    \lambda^{*(i)}_r.{\bar\lambda}^{*(k-i)}_{t-r}. N_i. z_i &
      \text{if } r \leq s_i, \; t\leq s_{k-i}, \\
      \lambda^{*(i)}_r. {\bar\lambda}^{(k-i)}_{t-r}.z_i  &
      \text{if } r \leq s_i, \; t > s_{k-i}, \\
      \lambda^{(i)}_r. {\bar\lambda}^{*(k-i)}_{t-r}. z_i &
      \text{if } r > s_i, \; t \leq s_{k-i}. 
  \end{array} \right.
$$
Here we have

\[ \begin{array}{lr}
 \Lambda^{*(3,5)}_{0,5} = \lambda^{*(3)}_0. \bar\lambda^{(5)}_5. z_3
                   = 57z_3, &
\Lambda^{*(3,5)}_{1,4} =\lambda^{*(3)}_1. \bar\lambda^{(5)}_4. z_3
                   = 9\times15z_3, \\ 
\Lambda^{*(3,5)}_{2,3} = \lambda^{*(3)}_2. \bar\lambda^{(5)}_3. z_3
                       =  120z_3, & 
\Lambda^{*(3,5)}_{3,2} = \lambda^{(3)}_3. \bar\lambda^{*(5)}_2. z_3
                    =  40z_3,\\
\Lambda^{*(3,5)}_{4,1} = \Lambda^{*(3,5)}_{5,0} = 0. &
\end{array}   \]

\[ \begin{array}{lr}
 \Lambda^{*(5,3)}_{0,5} = \Lambda^{*(5,3)}_{1,4}  =  0,  &
\Lambda^{*(5,3)}_{2,3} = \lambda^{*(5)}_2. \bar\lambda^{(3)}_3. z_5
                   = 40z_5, \\
\Lambda^{*(5,3)}_{3,2} = \lambda^{(5)}_3. \bar\lambda^{*(3)}_2. z_5
                    =  120z_5, &
\Lambda^{*(5,3)}_{4,1} = \lambda^{(5)}_4. \bar\lambda^{*(3)}_1. z_5
                      =  15\times 9 z_5, \\ 
\Lambda^{*(5,3)}_{5,0} = \lambda^{(5)}_5. \bar\lambda^{*(3)}_0. z_5
                       =  57z_5. & 
\end{array}  \]

\[ \begin{array}{lr}
 \Lambda^{*(4,4)}_{0,5} = \Lambda^{*(4,4)}_{5,0} =  0,  &
\Lambda^{*(4,4)}_{1,4} = \lambda^{*(4)}_1. \bar\lambda^{(4)}_4. z_4
                       = 48 z_4, \\
\Lambda^{*(4,4)}_{2,3} = \lambda^{*(4)}_2. \bar\lambda^{(4)}_3. z_4
                       =  8 \times 16 z_4, &
\Lambda^{*(4,4)}_{3,2} = \lambda^{(4)}_3. \bar\lambda^{*(4)}_2. z_4
                       =  16 \times 8 z_4, \\
\Lambda^{*(4,4)}_{4,1} = \lambda^{(4)}_4. \bar\lambda^{*(4)}_1. z_4
                       =  48 z_4. &
\end{array}   \]

It follows that
\begin{eqnarray*}
L_{0,5} &=&  u_0\bar{\lambda}^{(8)}_5 + u_1 19\bar{\lambda}^{(7)}_5 
            + u_2 171\bar{\lambda}^{(6)}_5 + u_3 57z_3, \\ 
L_{1,4} &=& u_1 5\bar{\lambda}^{(7)}_5 
            + u_2 9 \times 15 \bar{\lambda}^{(6)}_5 +
            u_3 9\times 15 z_3 + u_4 48z4 ,\\
L_{2,3} &=& u_2 40\bar{\lambda}^{(6)}_5 + u_3 120z_3
           + u_4 8\times 16 z_4 + u_5 40z5, \\
L_{3,2} &=&  u_6 40\lambda^{(6)}_5 + u_5 120 z_5
            + u_4 16\times 8z4 + u_3 40z_3, \\
L_{4,1} &=&  u_7 5\lambda^{(7)}_5  + u_6 15 \times 9 \lambda^{(6)}_5
       + u_5  15\times 9 z_5 + u_4 48z4, \\
L_{5,0} &=&  u_8 \lambda^{(8)}_5 + u_7 19\lambda^{(7)}_5  
        + u_6 171\lambda^{(6)}_5  + u_5 57z_5. 
 \end{eqnarray*}
Each set of values of $u_i\in \{0,1\}$, $i=0, \ldots, 8$;
 $z3, z4, z5 =1, \ldots, 17$; $\lambda^{(j)}_5$ and  
$\bar{\lambda}^{(j)}_5$, $j=6,7,8$ for which the equalities
$$L_{0,5} = L_{1,4} = L_{2,3} = L_{3,2} = L_{4,1} = L_{5,0} := \Lambda$$
is satisfied for a positive integer $\Lambda$ will yield
a simple $5-(38,8,\Lambda)$ design. Recall that a $5-(38,8,\Lambda)$
can be written as $5-(38,8,m4)$ with $\lambda_{\text{min}}=4$ and
$\lambda_{\text{max}}=5456$. 
Thus $\text{LIM}=\lfloor 5456/2*4 \rfloor = 682.$ 
By solving the equalities above we obtain for all $m4 \leq 1364$. 
Altogether 33 values for $m$ have been found, of which 16 values of 
$m \leq \text{LIM}$. 
Since, not all simple $5-(19,i,\lambda^{(i)}_5)$
designs are known to exist, for example, $5-(19,7, m7)$ designs
are known for $m= 4,5,6,7,8,9,13$ only, we just obtain the following
5 new simple $5-(38,8, m4)$ designs for $m=280, 488, 524, 560, 560$
(the number 560 repeats twice, as we have two distinct
 non isomorphic solutions for this value of $m$). The details of these
 5 constructed designs are given in Table 1.

\begin{table}[h!]
\begin{center}
\caption {Constructed simple $5-(38,8,\Lambda)$ designs} \label{Table1}
   \begin{tabular}{|l||rrrrrr|} \hline
    $m$ & $z_3$ & $z_4$ & $z_5$ & 
       $\lambda^{(6)}_5$ & $\lambda^{(7)}_5$ & $\lambda^{(8)}_5$ \\ 
        \hline \hline  
    280  &  7  &   0    &   7   &     0  &     35   &     56 \\
    488  &  8  &   4    &   8   &     4  &     28   &    280  \\
    524  &  6  &   7    &   6   &     6  &     28   &    196 \\
    560  &  4  &  10    &   4   &     8  &     28   &    112 \\
    560  &  9  &   5    &   9   &     4  &     49   &    112  \\ \hline
   \end{tabular}
  \end{center} 
 \end{table}
   An entry $0$ in a column of the table
   implies that $u_i=0$, otherwise $u_i=1$. 
   Here we have $\lambda^{(j)}_5 = \bar{\lambda}^{(j)}_5$, $j=6,7,8$
   for all these solutions.

\subsubsection{ Simple $5-(38,k,\Lambda)$ designs with $k=9, 10$}   
 Again we assume that $v_1=v_2=19$ for the construction of
 simple $5-(38,k,\Lambda)$ designs with $k=9, 10.$

 \begin{itemize}
  
  \item  For construction of $5-(38,9, \Lambda)= 5-(38,9, m30)$ designs 
        with $\mbox{LIM}= 682$, we make use of the large sets
        $\text{LS}[17](2,i,19)$, $i=3,4,5,6$, i.e. the 2-resolution
        of the complete designs $i-(19,i,1)$ with
        resolution class number $N_i=17$. Thus, we have 
        $R=\{(3,6), (4,5), (5,4), (6,3) \}.$ And the equalities 
        $L_{r,t-r}$ are the following.
 \begin{eqnarray*}
 L_{0,5} & = & u_0\bar\lambda^{(9)}_5 + u_1 19\bar\lambda^{(8)}_5 
         +u_2 171\bar\lambda^{(7)}_5  + u_3 57\times 14 z_3 + u_4  228 z_4,\\    L_{1,4} & = & u_1 15\bar\lambda^{(8)}_5/4 + u_2 18\times 5\bar\lambda^{(7)}_5 
         + u_3 9\times 105 z_3 + u_4 48\times 15 z_4 + u_5 180 z_5, \\
 L_{2,3} & = & u_2 20\bar\lambda^{(7)}_5 + u_3 560 z_3 + u_4 8\times 120 z_4
                +u_5 40\times 16 z_5 + u_6 140 z_6, \\
 L_{3,2} & = & u_7 20\lambda^{(7)}_5 + u_6 560 z_6 +u_5 120\times 8 z_5
                +u_4 16\times 40 z_4 + u_3 140 z_3, \\ 
 L_{4,1} & = & u_8 15\lambda^{(8)}_5/4 
                +u_7  5\times 18\lambda^{(7)}_5
                +u_6 105\times 9 z_6 + u_5 15\times 48 z_5 + u_4 180 z_4, \\
 L_{5,0} & = & u_9 \lambda^{(9)}_5 + u_8 19\lambda^{(8)}_5 
                +u_7 171\lambda^{(7)}_5 + u_6 14\times 57 z_6 
                + u_5 228 z_5.
 \end{eqnarray*}
      Solving the equalities 
      $L_{0,5}= L_{1,4}= L_{2,3}= L_{3,2}= L_{4,1}= L_{5,0} = \Lambda$
      for $\Lambda > 0$  with respect to $z_i=1, \ldots, 17$ 
      we obtain 20 values for 
      $m$ with $m\leq \mbox{LIM}$ leading to simple 
      $5-(38,9, \Lambda)=5-(38,9, m30)$ designs.
      Of which 14 designs can be constructed whose details
      are given in Table 2.

\begin{table}[h!]
\caption{Constructed simple $5-(38,9,\Lambda)$ designs}\label{Table2}             \begin{center}
   \begin{tabular}{|l||rrrr rrr rrr|} \hline
    $m$ & $z_3$ & $z_4$ & $z_5$ & $z_6$ &
       $\lambda^{(7)}_5$ & $\lambda^{(8)}_5$ & $\lambda^{(9)}_5$ &
    $\bar\lambda^{(7)}_5$ & $\bar\lambda^{(8)}_5$ & $\bar\lambda^{(9)}_5$ \\ 
        \hline \hline 
 100 & 2  &  1 &  1  &  2  &  0  &   56 &   112  &   0 &   56 &  112 \\
 200 & 4  &  2 &  2  &  4  &  0  &  112 &   224  &   0 &  112 &  224 \\
 300 & 6  &  3 &  3  &  6  &  0  &  168 &   336  &   0 &  168 &  336 \\
 400 & 8  &  4 &  4  &  8  &  0  &  224 &   448  &   0 &  224 &  448 \\
 402 & 5  &  5 &  5  &  5  & 28  &   84 &   546  &  28 &   84 &  546 \\
 500 & 10 &  5 &  5  & 10  &  0  &  280 &   560  &   0 &  280 &  560 \\
 502 & 7  &  6 &  6  &  7  & 28  &  140 &   658  &  28 &  140 &  658 \\
 504 & 4  &  7 &  7  &  4  & 56  &    0 &   756  &  56 &    0 &  756 \\
 582 & 10  & 4 & 11  &  3  &  28 &  168 &  588   &  63 &   84 &  189 \\
 602 & 9  &  7 &  7  &  9  & 28  &  196 &   770  &  28 &  196 &  770 \\ 
 604 & 6  &  8 &  8  &  6  & 56  &   56 &   868  &  56 &   56 &  868 \\
 660 & 9  &  8 &  8  &  9  & 35  &  252 &    21  &  35 &  252 &   21 \\
 680 & 8 &  11 &  4  & 15  &  0  &  364 &   602  &  35 & 280  &  203 \\   
 682 & 5  & 12 &  5  & 12  & 28  &  224 &   700  &  63 & 140  &  301 \\ \hline 
  \end{tabular}
  \end{center} 
\end{table}
   
It should be noted that when applying the basic construction 
for $t=5$, $v_1=v_2=19$ and $k=8,9$ we only obtain the trivial solutions,
namely the complete $5-(38,8,1364\times 4)$
and $5-(38,9,1364\times 30)$ designs. This could be explained as follows.

In general, if $k \leq 2t-1$, then one of the designs in each pair
$(D_i, \bar D_{k-i})$ is either the empty or the trivial design
and at least one pair having both the trivial designs, therefore 
it leaves little room for the basic construction to produce a 
non-trivial solution, unless many pairs are unused, i.e. 
$u_i=0.$  
The construction using resolution indeed
makes more room to create non-trivial solutions, as we have seen in the 
above examples.

  \item  For construction of $5-(38,10, \Lambda)= 5-(38,10, m6)$ designs 
        with $\mbox{LIM}= 19778$, we again employ the 2-resolution
        of the complete designs $i-(19,i,1)$ for $i=3,4,5,6,7$ with
        resolution class number $N_i=17$. Here, 
        $$R=\{(3,7), (4,6), (5,5), (6,4), (7,3) \}.$$  And we have

 \begin{eqnarray*}
  L_{0,5} & = & u_0 \bar\lambda^{(10)}_5 + u_1 19\bar\lambda^{(9)}_5 
               + u_2 171\bar\lambda^{(8)}_5  + u_3 57 \times 91 z_3  
               + u_4 228 \times 14 z_4 + u_5  684 z_5, \\
  L_{1,4} & = & u_1 3\bar\lambda^{(9)}_5 + u_2 18 \times 
               15\bar\lambda^{(8)}_5/4 + u_3 9 \times 455 z_3 
               +u_4 48 \times 105 z_4  + u_5 180 \times 15 z_5 \\ 
           &  & + u_6  504 z_6, \\
  L_{2,3} & = & u_2 12 \bar\lambda^{(8)}_5 + u_3 1820 z_3
               + u_4 8 \times 560 z_4  +u_5 40 \times 120 z_5
               + u_6 140 \times 16 z_6  \\ 
          &  & + u_7 364 z_7, \\
  L_{3,2} & = & u_8 12 \lambda^{(8)}_5 + u_7 1820 z_7
               + u_6 560 \times 8 z_6  +u_5 120 \times 40 z_5
               + u_4 16 \times 140 z_4 \\
           &  & + u_3 364 z_3, \\
  L_{4,1} & = & u_9 3 \lambda^{(9)}_5 + u_8 15 \times 18\lambda^{(8)}_5/4
               + u_7 455 \times 9 z_7 +u_6 105 \times 48 z_6 
               + u_5 15 \times 180 z_5 \\
           &  & + u_4  504 z_4, \\
  L_{5,0} & = & u_{10} \lambda^{(10)}_5 + u_9 19\lambda^{(9)}_5 
               + u_8 171\lambda^{(8)}_5  +u_7 91 \times 57 z_7  
               + u_6 14 \times 228 z_6 + u_5  684 z_5.
 \end{eqnarray*}
 
      Solving the equalities 
      $L_{0,5}= L_{1,4}= L_{2,3}= L_{3,2}= L_{4,1}= L_{5,0} = \Lambda$
      for $\Lambda > 0$  with respect to $z_i=1, \ldots, 17$ 
      we obtain an entire number of 479 solutions, of which 
      239 have $m\leq \mbox{LIM}$. From these 239 parameters
      131 simple $5-(38,10, m6)$ designs have been shown to exist.
      The values of $m$ for these designs are

\[
\begin{array}{rrrrrrrr}
   12768 &  17416 &   2604 &   6076  &  7252 &  10724 &  13668 &  15108 \\
   15372 &  18580 &  18844 &   3768  &  6976 &   8416 &   8680 &  11624 \\
   11888 &  12152 &  16272 &  16536  & 16800 &  19744 &   4932 &   8404 \\ 
    9580 &   9844 &  12788 &  13052  & 13316 &  13580 &  17172 &  17436 \\
   17700 &  17964 &  18228 &   6096  &  9040 &  10480 &  10744 &  11008 \\
   13952 &  11536 &  14216 &  14480  & 15920 &  18600 &  18864 &  19128 \\
    7260 &  11644 &  11908 &  12172  & 14852 &  15116 &  15380 &  15644 \\
   16556 &  16820 &  19500 &  17084  & 19764 &   8424 &  11368 &  12544 \\
   12808 &  13072 &  13336 &  16016  & 16280 &  16544 &  16808 &  17984 \\
    9060 &   9588 &  13972 &  14236  & 16916 &  14500 &  17180 &  17444 \\
   17708 &  18884 &  19148 &  10224  & 10752 &  14872 &  15136 &  15400 \\
   18080 &  15664 &  18344 &  18608  & 19520 &  11388 &  11916 &  16036 \\
   16300 &  16564 &  19244 &  16828  & 19508 &  19772 &  12552 &  13080 \\
   17200 &  17464 &  17728 &  13716  & 14244 &  18100 &  18364 &  18628 \\
   18892 &  19156 &  14880 &  15408  & 19264 &  19528 &  16044 &  17208 \\
   18384 &  18372 &  19536 &  16844  & 11316 &  13908 &  14280 &  14808 \\
   19720 &  16872 &  17772 &    &  &  &  &
\end{array}    \]

\end{itemize}
    
 Here are two examples: 
 \begin{itemize}
  \item
  $5-(38,10,2604\times6)$ with $z_3=1$, $z_4=2$, $z_6=2$,
   $z_7=1$, $\bar\lambda^{(9)}_5=147$, $\bar\lambda^{(10)}_5=1260$,
   $u_2=u_5=u_8=0$ and $\lambda^{(i)}_5 =\bar\lambda^{(i)}_5$
    for $i=9,10.$ 
  \item
  $5-(38,10,11316\times6)$ with $z_3=2$, $z_4=8$, $z_5=2$, $z_6=7$,
   $z_7=4$, $\bar\lambda^{(8)}_5=140$, $\bar\lambda^{(9)}_5=336$, 
   $\bar\lambda^{(10)}_5=294$,
   $\lambda^{(8)}_5 = 84$, $\lambda^{(9)}_5=378$, $\lambda^{(10)}_5=1890.$

 \end{itemize}

On the other hand, when the basic construction is applied for
      this case (i.e. $v_1=v_2=19$ and $k=10$), we just obtain
      5 solutions with $m \leq \text{LIM}$.
\begin{rem}
{\rm  
  \begin{enumerate}
 \item
  It should be noted that when $v_1=v_2$, any solution with
  $\lambda^{(i)}_t \not= \bar{\lambda}^{(i)}_t$ 
  will appear twice by reason of symmetry, since  $\lambda^{(i)}_t$
  and $\bar{\lambda}^{(i)}_t$ may be interchanged. These two solutions
  are indeed the same. This fact should
  be taken into account by counting the number of solutions.

\item
 Up to now the number of known simple designs 
 for $5-(38,k,\Lambda)$ with $k=8,9,10$ are 
 8, 14, and 23 respectively, see \cite{handbook07}, for instance.
 For $k=8,9$ all the parameters of the constructed designs differ
 from the known ones. For $k=10$, only one of the 23 parameters
 of the known designs does appear in the list of 131 constructed
 designs, namely the parameters $5-(38,10, 11368\times 6)$. However,
 it is not known whether the corresponding designs are isomorphic. 

\end{enumerate}
}
\end{rem}

\subsection{Some further results of applications}
We briefly record some further examples of simple $t$-designs
for $t=4,5,6$ by using Theorem \ref{main}. 

\subsubsection{$t=4$}
Following are several small parameters for $t=4$.
\begin{enumerate}
\item
 $4-(26,8,m35)$:  Using $v_1=v_2=13$ and a subset
 $R= \{(3,5), (4,4), (5,3)\}$ of pairs of designs having $s_i$-resolution
 derived from  $\text{LS}[55](2,4,13)$ and 
 $\text{LS}[11](2,i,13)$ for $i=3,5$. 
 There are 3 non-trivial solutions of Eq(\ref{Equalities2}) 
 with $m=44,66$ satisfying  $m \leq \text{LIM}(=104).$  
 A design with $m=44$ is known.
 The two solutions for $m=66$ are non-isomorphic and new. These are
 \begin{itemize}
\item $u_4=0$, $z_3=z_5=7$, $\lambda^{(7)}_4=42$, $\lambda^{(8)}_4=126$,
   $u_2=u_6=0$, and $\bar\lambda^{(i)}_4= \lambda^{(i)}_4$ for $i=7,8$.

 \item $z_4=24$, $z_3=z_5=2$, $\lambda^{(6)}_4=18$, $\lambda^{(8)}_4=126$,
   $u_1=u_7=0$, and $\bar\lambda^{(i)}_4= \lambda^{(i)}_4$ for $i=6,8$.
 \end{itemize} 
 The basic construction for $4-(26,8,m35)$ with $v_1=v_2=13$
 only yields the trivial solution.

\item 
 $4-(28,9,m168)$:  Using $v_1=v_2=14$ and a subset of resolution pairs
 $R= \{(4,5),(5,4)\}$
 derived from  $\text{LS}[11](2,i,14)$ for $i=4,5$.
 There is a unique non-trivial solution of Eq(\ref{Equalities2})
 with  $m=110$ satisfying $m \leq \text{LIM}(=126).$ 
 This solution with
 $z_4=z_5=4$, $u_2=u_7=0$, $\lambda^{(6)}_4=30$, $\lambda^{(8)}_4=210$, 
  $\lambda^{(9)}_4=252$,
 and $\bar\lambda^{(i)}_4= \lambda^{(i)}_4$ for $i=6,8,9$ 
 yields a new design. 

\item
 $4-(30,7,m20)$:  Using $v_1=v_2=15$ and a subset of resolution pairs
 $R= \{(3,4),(4,3)\}$ 
 derived from  $\text{LS}[13](2,i,15)$ for $i=3,4$.
 There are 3 non-trivial solutions of Eq(\ref{Equalities2})
 with $m=39,52,65$ satisfying $m \leq \text{LIM}(=65).$  
 The solution for $m=52$ with
 $z_3=z_5=5$, $\lambda^{(5)}_4=5$, $\lambda^{(6)}_4=15$, 
 $\lambda^{(7)}_4=115$,
 and $\bar\lambda^{(i)}_4= \lambda^{(i)}_4$ for $i=5,6,7$
 gives a new design.

\end{enumerate}

\subsubsection{$t=5$}
\begin{enumerate}
\item
$5-(36,10,m63)$:  Using $v_1=v_2=18$ and a subset of resolution pairs
 $R= \{(5,5)\}$
 derived from  $\text{LS}[7](2,5,18)$. There are 164 non-trivial solutions of 
 Eq(\ref{Equalities2})
 with  $m \leq \text{LIM}(=1348).$  
 Of which 37 are shown to exist.
 It is interesting to remark that these 37 designs
 include the 10 designs constructed using the  basic construction
 \cite{TvT2016}. Actually, 27 new designs with
 parameters $5-(36,10,m63)$ have been obtained. These are
 \begin{eqnarray*} 
 m &= & 611, 818, 921, 945, 969, 1048, 1072, 911, 934, 1094,
  1197, 1221, 1245, 1269,\\
  & & 1324, 1325, 1348, 1187, 1210, 1234, 1337, 1152, 1176, 1200,
  1224, 1303, 1131.  
\end{eqnarray*}

\item 
$5-(37,8,m40)$:  Using $v_1=13$, $v_2=24$ and a subset of resolution pairs
 $R= \{(3,5),(4,4), (5,3) \}$
 derived from  $\text{LS}[11](2,i,13)$, $\text{LS}[11](2,i,24)$ 
 for $i=3,4,5$. 
 There is a unique non-trivial solution of
 Eq(\ref{Equalities2}) with $m=55$ such that
 $m \leq \text{LIM}(=62).$  This solution
 with $z_3=2$,$z_5=2$, $z_4=8$, 
 $\lambda^{(6)}_5=4$, $\lambda^{(7)}_5=28$, $\lambda^{(8)}_5=56$,
 $\bar\lambda^{(6)}_5= 13$, $\bar\lambda^{(7)}_5= 36$,
 $\bar\lambda^{(7)}_5= 666$
 gives a new design.

\item
$5-(37,9,m10)$:  Using $v_1=13$, $v_2=24$ and a subset of resolution pairs
 $R= \{(3,6),(4,5), (5,4),(6,3) \}$
 derived from  $\text{LS}[11](2,i,13)$, 
 $\text{LS}[11](2,i,24)$ for $i=3,4,5,6$.
 There is a unique non-trivial solution of
 Eq(\ref{Equalities2}) with $m=874$ such that
 $m \leq \text{LIM}(=1798).$  This solution
 with $z_3=2$, $z_4=2$, $z_5=4$, $z_6=1$,
 $\lambda^{(7)}_5=14$, $u_8=u_9=0$, 
 $\bar\lambda^{(6)}_5= 72$, $\bar\lambda^{(7)}_5= 30$,
 $\bar\lambda^{(8)}_5= 1980$
 gives a new design.
 
\item 
$5-(44,8,m)$:  Using $v_1=v_2=22$ and a subset of resolution pairs
 $R= \{(4,4)\}$
 derived from  $\text{LS}[19](2,5,22)$. There are 9 non-trivial solutions of
 Eq(\ref{Equalities2})
 with  $m \leq \text{LIM}(=4569).$ 
 Of which one design with $m=3344$ and
 $u_3=u_5=0$, $z_4=4$,
 $\lambda^{(6)}_5=12$, $\lambda^{(7)}_5=16$, $\lambda^{(8)}_5=220$,
 and $\bar\lambda^{(i)}_5= \lambda^{(i)}_5$ for $i=6,7,8,$
 is shown to exist. 
\item 
$5-(46,10,m2)$:  Using $v_1=v_2=23$ and a subset of resolution pairs
 $R= \{(4,6), (5,5), (6,4)\}$
 derived from  $\text{LS}[133](2,5,23)$ and 
  $\text{LS}[7](2,i,23)$ for $i=4,6$. 
There are 3986 non-trivial solutions of
 Eq(\ref{Equalities2})
 with  $m \leq \text{LIM}(=187349).$ 
 Of which 176 designs are shown to exist with the following values
 of $m$.

\[
\begin{array}{rrrrrrrr}
   65246  &   75487  &   73758  &   83999  &   86526  &   94240  &   96140  &   96767 \\
  106381  &  107008  &  116622  &  117021  &  125134  &  123405  &  127262  &  139004 \\
  137503  &  142633  &  139403  &  143887  &  149644  &  151772  &  159885  &  162013 \\
  164540  &  166668  &  174781  &  185497  &   59014  &   79667  &   78166  &   89908 \\
   88179  &   92435  &   98420  &  102676  &  108661  &  110561  &  120802  &  126160 \\
  122930  &  125058  &  131043  &  136401  &  133171  &  139555  &  137826  &  145540 \\
  153425  &  157054  &  158308  &  175807  &  174705  &  184319  &  182818  &  184946 \\
   77064  &   94088  &   99446  &   96216  &  104329  &  109687  &  102600  &  119928 \\
  112841  &  116698  &  121828  &  115368  &  123082  &  132069  &  124982  &  125609 \\
  130967  &  142310  &  135223  &  140581  &  135850  &  141208  &  145464  &  150822 \\
  145863  &  153976  &  156104  &  163590  &  167846  &  166345  &  171475  &  168245 \\
  178486  &  183844  &  180614  &  185972  &   86526  &  105355  &  113867  &  112366 \\
  118750  &  124108  &  126635  &  134349  &  127262  &  136249  &  134748  &  131518 \\
  136876  &  137503  &  146490  &  139403  &  147117  &  140030  &  156731  &  149644 \\
  155002  &  151772  &  157130  &  150271  &  153900  &  159885  &  165243  &  162013 \\
  167371  &  166668  &  177612  &  174382  &  179512  &  182267  &  185896  &  187150 \\
  121505  &  119776  &  128288  &  132544  &  138529  &  142785  &  148770  &  145540 \\
  150670  &  160911  &  163039  &  159809  &  171152  &  169423  &  164692  &  170050 \\
  168321  &  179664  &  174705  &  182818  &  184946  &  135850  &  148694  &  152950 \\
  158707  &  163191  &  165091  &  163590  &  165718  &  171076  &  175332  &  179189 \\
  175959  &  181317  &  174230  &  185573  &  178486  &  183844  &  180614  &  185972 \\ 
  184471  &  182742  &  176985  &  164540  &  177612  &  179512  &  187226  &  173052  
\end{array}
 \]

 Here is an
 example with $m=59014$: $z_4=z_6=1$, $z_5=20$, $u_2=u_8=0$,
 $\lambda^{(7)}_5=36$, $\lambda^{(9)}_5=810$, $\lambda^{(10)}_5=7812$,
 and $\bar\lambda^{(i)}_5= \lambda^{(i)}_5$ for $i=7,9,10.$

\end{enumerate}

\subsubsection{$t=6$}
Following are some examples for $t=6$.
\begin{enumerate}
\item
 $6-(38,10,m10)$:  Using $v_1=v_2=19$ and a subset of resolution pairs
 $R= \{(4,6),(5,5),(6,4)\}$
 derived from  $\text{LS}[4](3,i,19)$ for $i=4,5,6$. 
 There are 4 non-trivial solutions of 
 Eq(\ref{Equalities2}) with $m= 1360, 892, 1340, 1788$ for
$m \leq \text{LIM}(=1798).$

\item
 $6-(46,12,m420)$:  Using $v_1=v_2=23$ and a subset of resolution pairs
 $R= \{(6,6)\}$ derived from  $\text{LS}[3](3,6,23)$.                     
 There are 2 non-trivial solutions of  
 Eq(\ref{Equalities2}) with $m= 3363, 3819$ for 
 $m \leq \text{LIM}(=4569).$  The solution for $m=3363$
 has $z_6=1$, $\lambda^{(7)}_5=7$, $\lambda^{(8)}_5=40$, $\lambda^{(9)}_5=340$,
 $\lambda^{(10)}_5=350$, $\lambda^{(11)}_5=4046$, $\lambda^{(12)}_5=5320$,
 and $\bar\lambda^{(i)}_5= \lambda^{(i)}_5$ for $i=7,8,9,10,11,12.$
 All the ingredient designs corresponding to $m=3363$ exist except that
 the existence of a $6-(23,10,5\times 70)$ design is still in doubt. 
 So, we would have a $6-(46,12, 3364\times 420)$ design 
 if a $6-(23,10,5\times 70)$ design would exist.

\item
 $6-(50,12,m308)$:  Using $v_1=v_2=25$ and a subset of resolution pairs
 $R= \{(6,6)\}$ derived from  $\text{LS}[7](3,6,25)$.
 There are 195 non-trivial solutions of
 Eq(\ref{Equalities2}) for $m \leq \text{LIM}(=11459).$

\end{enumerate}

\section{Conclusion}
We have presented a recursive construction for simple $t-$designs
by using the concept of resolution. This may be viewed as an extension
of the basic construction as shown in our previous paper. 
The $s$-resolutions of trivial $t$-designs are equivalent to 
the large sets of $s$-designs, which have been extensively studied.
Since our construction does not exclude the use of trivial designs as
ingredients, we have restricted its applications to resolutions of the
trivial ingredient designs only. In spite of this fact, the construction
still produces a large number of new simple $t$-designs.
We do not know any $s$-resolutions of non-trivial 
$t$-designs for $t \geq 4$ and $s \geq 2$. However,
we strongly believe that the construction would unfold its full impact
when we would gain more knowledge about resolutions of 
non-trivial $t$-designs.

\end{document}